\documentclass{amsart}
\usepackage{amsmath,amsthm,amssymb,amscd}
\theoremstyle{plain}

\newfont{\cyrrlarge}{wncyr10}
\def\Sha{\mbox{\cyrrlarge Sh}}

\newfont{\cyrrmed}{wncyr8}
\def\Shamed{\mbox{\cyrrmed Sh}}

\newfont{\cyrr}{wncyr5}

\newtheorem{thm}{Theorem}[section]
\newtheorem{prop}[thm]{Proposition}
\newtheorem{lemma}[thm]{Lemma}

\theoremstyle{definition}

\newtheorem{defn}[thm]{Definition}
\newtheorem{remark}[thm]{Remark}

\newcommand{\Z}{{\bf Z}}
\newcommand{\Q}{{\bf Q}}

\newcommand{\Zp}{\Z_p}
\newcommand{\Zpx}{\Z_p^\times}
\newcommand{\Qp}{\Q_p}

\newcommand{\Fp}{{\bf F}_p}

\newcommand{\Qinf}{\Q_\infty}
\newcommand{\Qinfp}{\Q_{\infty,p}}
\newcommand{\Qn}{\Q_n}

\newcommand{\surj}{\twoheadrightarrow}
\newcommand{\maps}{\longrightarrow}

\newcommand{\ve}{\varepsilon}

\newcommand{\F}{\mathcal F}

\newcommand{\pin}{\pi_{n/n-1}}
\newcommand{\nun}{\nu_{n-1/n}}

\DeclareMathOperator{\ord}{ord}
\newcommand{\ordp}{\ord_p}

\DeclareMathOperator{\Tam}{Tam}
\DeclareMathOperator{\Gal}{Gal}
\DeclareMathOperator{\Tr}{Tr}
\DeclareMathOperator{\Sel}{Sel}

\DeclareMathOperator{\im}{im}

\newcommand{\dual}{\wedge}

\title{An Algebraic Version of a Theorem of Kurihara}

\author{Robert Pollack}
\address{Department of Mathematics,
Boston University, Boston MA 02215, USA}
\email{rpollack@math.bu.edu}
\thanks{The author was supported by an NSF Postdoctoral Fellowship.}
\dedicatory{In memory of Professor Arnold Ross}

\begin{document}

\begin{abstract}
Let $E/\Q$ be an elliptic curve and let $p$ be an odd supersingular
prime for $E$. 
In this article, we study the simplest case of Iwasawa theory for
elliptic curves, namely when $E(\Q)$ is finite, $\Shamed(E/\Q)$ has no $p$-torsion
and the Tamagawa factors for $E$ are all prime to $p$. 
Under these hypotheses, we prove that $E(\Qn)$ is finite and make 
precise statements about the size and structure of the $p$-power part
of $\Shamed(E/\Qn)$.  Here
$\Qn$ is the $n$-th step in the cyclotomic $\Zp$-extension of
$\Q$.  
\end{abstract}

\maketitle

\section{Introduction}
\label{section:intro}

Let $E/\Q$ be an elliptic curve with good supersingular reduction at an odd prime $p$.
Let $\Qinf$ be the cyclotomic $\Zp$-extension of $\Q$ with subfields $\Qn$ of degree
$p^n$.  In \cite{Kurihara}, Kurihara proved precise statements about the size and
the structure of the $p$-part of the Tate-Shafarevich group $\Sha(E/\Qn)$
when $\ordp(L(E,1)/\Omega_E) = 0$ and when the Galois representation on the $p$-torsion
is surjective.  His proof made deep use of Kato's Euler system for the Tate
module of $E$ (and hence the need for an assumption on the Galois representation).

In this paper, we offer a completely algebraic proof of a variant of a theorem of Kurihara
(see \cite[Theorem 0.1]{Kurihara}) 
where his analytic assumptions are converted to algebraic ones (equivalent under
the Birch and Swinnerton-Dyer conjecture).  
Before stating the result, we fix some notation.
Set $\Gamma = \Gal(\Qinf/\Q)$, $\Gamma_n = \Gal(\Qinf/\Qn)$ and $G_n = \Gal(\Qn/\Q)$.  
Let $\Lambda_n = \Zp[G_n]$ be the group algebra at level $n$ and
$\Lambda = \Zp[[\Gamma]]$ be the Iwasawa algebra.
For a $\Zp$-module $M$, denote by $M^\dual$ its Pontrjagin dual.

\begin{thm}
\label{thm:main}
Let $E/\Q$ be an elliptic curve with $p$ an odd prime of good supersingular reduction.
Assume that 
\begin{enumerate}

\item $E(\Q)$ is finite
\item $p \nmid \Tam(E/\Q)$
\item $\Sha(E/\Q)[p]=0$.

\end{enumerate}
Then 
\begin{enumerate}

\item $E(\Qn)$ is finite for all $n \geq 0$.
\item $\ordp(\#\Sha(E/\Qn)) = e_n$ where $e_0=e_1=0$ and
$$
e_n = 
\begin{cases}
p^{n-1} + p^{n-3} + \cdots + p - \frac{n}{2} \text{~for~even~} n \geq 2 \\
p^{n-1} + p^{n-3} + \cdots + p^2 - \frac{n-1}{2} \text{~for~odd~} n \geq 3.
\end{cases}
$$
\item \label{item:struct}
When $a_p=0$, we have
$$
\Sha(E/\Qn)[p^\infty]^\dual \cong \Lambda_n / (J_n^+ + J_n^-)
$$
as $\Zp[G_n]$-modules where
\begin{align*}
J_n^\pm := \{ &f \in \Lambda_n : \chi(f) = 0 \text{~for~} \chi \text{~a~char.~of~}
G_n \text{~of~even~(resp.~odd)~order} \}.
\end{align*}

\end{enumerate}

\end{thm}

\begin{remark}
The above theorem is false for $p=2$.  If $E=X_0(19)$ then $E(\Q)$ is finite,
$\Tam(E/\Q)$ is odd
and $\Sha(E/\Q)[2]=0$.  However, $E(\Q(\sqrt{2}))$ is infinite and $\Q(\sqrt{2})$
is the first step in the cyclotomic $\Z_2$-extension.
\end{remark}

\begin{remark}
The conclusion of Theorem \ref{thm:main} is identical to Kurihara's theorem; it
is only the hypotheses that have changed.
For supersingular $p$, 
the Birch and Swinnerton-Dyer conjecture predicts that
$\ordp\left( L(E,1) / \Omega_E  \right) = 0$ if and only if
$E(\Q)$ is finite, $p \nmid \Tam(E/\Q)$ and $\Sha(E/\Q)[p]=0$.
The ``if part'' is still unknown, but the ``only if'' part is known via Kato's Euler system
when the Galois representation on the $p$-torsion
is surjective.  Hence the above hypotheses are logically weaker than Kurihara's since
we make no assumptions on the Galois representation.  In particular, our results apply
to CM curves.
\end{remark}

The analogue of Theorem \ref{thm:main} in the ordinary case follows from
Mazur's control theorem. However, in the supersingular case the control
theorem fails (due to the triviality of the universal norms of the formal
group $\hat{E}/\Qp$ along the local cyclotomic $\Zp$-extension). We will
make a careful study of the how the control theorem fails in terms of
$\hat{E}$ and combining this with a precise enough description of this
formal group, we will be able to prove Theorem \ref{thm:main}.

These techniques are not new as they form the basis of Perrin-Riou's
construction of an algebraic $p$-adic $L$-function in \cite{PR:90}.  
Also, many of the calculations in this paper were inspired by the 
beautiful ideas of Kurihara in \cite{Kurihara}.
It should also be mentioned that similar results were announced
by Nasybullin \cite{Nasybullin} over twenty five years ago, but 
in his short paper no proofs were given. 

One advantage to the algebraic approach of this paper is that it can be generalized
more easily to $\Zp$-extensions of a number field that are not necessarily cyclotomic.
To successfully carry out such a generalization, the key local input that is needed is a 
good understanding of the Galois module structure of $\hat{E}$ along the $\Zp$-extensions 
of some finite extension of $\Qp$.  In a forthcoming paper with Adrian Iovita
(see \cite{Iovita}) a strong enough local result is obtained to generalize the results
of this paper to any $\Zp$-extension of a number field in which $p$ splits completely.

The format of the paper will be as follows: in the following section we
will implement the needed Iwasawa theory to precisely describe the failure
of the control theorem in terms of $\hat{E}$.  The third section will
state results of Kobayashi on the structure of $\hat{E}$ as a Galois
module.  In the fourth section, we will define $\mu$ and $\lambda$-invariants
of elements of $\Lambda_n$ and discuss their basic properties.
In the final section, we will perform the needed computations to
complete the proof of Theorem \ref{thm:main}.

~\\
\noindent
{\it Acknowledgments}: We would like to thank Ralph Greenberg and Adrian Iovita for
many interesting conversations on the Iwasawa theory of elliptic curves at supersingular
primes.  We also thank the anonymous referee for pointing out an
important detail that was not verified.

\section{Iwasawa theory}

Let $E/\Q$ be an elliptic curve, $p$ some prime of good reduction and $K$ some finite
extension of $\Q$.  We define the $p$-Selmer group of $E$ over $K$ by
$$
\Sel_p(E/K) = \ker \left( H^1(K,E[p^\infty]) \maps \prod_v H^1(K_v,E) \right)
$$
where $v$ runs over the places of $K$.  Also, define a looser Selmer group by dropping
the condition at $p$, i.e.
$$
\Sel_p'(E/K) = \ker \left( H^1(K,E[p^\infty]) \maps \prod_{v \nmid p} H^1(K_v,E) \right).
$$
We then have the following exact sequence relating these two Selmer groups:
\begin{eqnarray}
\label{eqn:map1}
0 \maps \Sel_p(E/K) \maps \Sel_p'(E/K) \maps \prod_{v \mid p} H^1(K_v,E)[p^\infty] .
\end{eqnarray}

For the infinite extension $\Qinf$ we define 
$\displaystyle \Sel_p(E/\Qinf) = \varinjlim \Sel_p(E/\Qn)$ and  
$\displaystyle \Sel_p'(E/\Qinf) = \varinjlim \Sel_p'(E/\Qn)$.  
As mentioned in the introduction, the control theorem for
$\Sel_p(E/\Qinf)$ fails for supersingular $p$.  However, the control theorem 
for $\Sel_p'(E/\Qinf)$ is always true.

\begin{thm}
\label{thm:control}
Let $p$ be a prime of good reduction for $E/\Q$.  Then the natural map
$$
\Sel_p'(E/\Qn) \maps \Sel_p'(E/\Qinf)^{\Gamma_n}
$$
has finite kernel and cokernel that are bounded independent of $n$.  

Moreover, if $E(\Q)[p]=0$, $p \nmid \Tam(E/\Q)$ and $a_p \not\equiv 1 \pmod{p}$
then the above map is an isomorphism.
\end{thm}

\begin{proof}
This theorem was originally proven by Mazur in \cite{Mazur}.  See also
\cite{Manin} and \cite[Chapter 3]{Greenberg:IEC} for an exposition of
this theorem that uses Galois cohomology instead of flat cohomology.
Note that in all of these papers the ordinary
hypothesis is only used in studying the primes over $p$.
Since we are dealing with $\Sel'$ and not $\Sel$ these proofs apply to our situation.
\end{proof}

We now work under the hypotheses of Theorem \ref{thm:main}, namely that
$p$ is supersingular for $E$, $E(\Q)$ is finite, $p \nmid \Tam(E/\Q)$ and
$\Sha(E/\Q)[p]=0$.  Since $p$ is supersingular, $a_p \not\equiv 1 \pmod{p}$ and
$E(\Q)[p]=0$.
Hence, the map in Theorem \ref{thm:control} is an isomorphism and 
(\ref{eqn:map1}) becomes
\begin{eqnarray}
\label{eqn:map2}
0 \maps \Sel_p(E/\Qn) \maps \Sel_p'(E/\Qinf)^{\Gamma_n} \maps H^1(\Q_{n,p},E)[p^\infty]
\end{eqnarray}
where $\Q_{n,p}$ denotes the completion of $\Qn$ at the unique prime over $p$.

The main reason for the failure of the control theorem in the supersingular
case is that the local condition defining the Selmer group at $p$ disappears
over $\Qinf$.  

\begin{prop}
\label{prop:deep}
For $p$ supersingular
$$
H^1(\Qinfp,E)[p^\infty] = 0
$$
and hence
$$
\Sel_p(E/\Qinf) = \Sel_p'(E/\Qinf).
$$
\end{prop}

\begin{proof}
By Tate local duality, the vanishing of $H^1(\Qinfp,E)[p^\infty]$ is equivalent to the
triviality of the universal norms of $\hat{E}$ along $\Qinfp/\Qp$.  This vanishing
of universal norms was originally proven by Hazewinkel in \cite{Hazewinkel}.
See \cite{C-G:Kummer} for a general discussion of this phenomenon for deeply
ramified extensions.
\end{proof}

Hence $X_\infty := \Sel_p(E/\Qinf)^\dual \cong \Sel'_p(E/\Qinf)^\dual$.
Dualizing (\ref{eqn:map2}) and applying Tate local duality yields
\begin{eqnarray}
\label{eqn:map3}
\hat{E}(\Q_{n,p}) \maps \left( X_\infty \right)_{\Gamma_n} \maps \Sel_p(E/\Qn)^\dual \maps 0
\end{eqnarray}
where $M_{\Gamma_n}$ denotes the $\Gamma_n$-coinvariants of $M$.  
The above sequence can be thought of as describing the failure of the control
theorem in terms of the formal group.

We make one last alteration of the above sequence by explicitly describing
$X_\infty$.  The following is well known, but we include a proof for completeness.

\begin{prop}
Under our hypotheses, $X_\infty$ is a free $\Lambda$-module of rank 1.
\end{prop}

\begin{proof} When $p$ is supersingular, it is always true that the
$\Lambda$-rank of $X_\infty$ is greater than or equal to 1 by a result of
Schneider (see \cite[Corollary 5]{Schneider}).  For a discussion of this
theorem using Galois cohomology rather than flat cohomology see
\cite[Proposition 2.6]{C-S:GCEC}.

Under our hypotheses, we prove an upper bound on the $\Lambda$-rank of $X_\infty$ and
establish that it is a free $\Lambda$-module.
Note that since $E(\Q)$ is finite and $\Sha(E/\Q)[p]=0$ we have that
$\Sel_p(E/\Q)=0$.  Hence, taking $n=0$ in (\ref{eqn:map3}) yields
$$
\hat{E}(\Qp) \surj \left( X_\infty \right)_{\Gamma}.
$$
Furthermore, $\hat{E}(\Qp) \cong \Zp$ and since $\left( X_\infty \right)_{\Gamma}$
is infinite the above map is an isomorphism.  A compact version of
Nakayama's lemma then implies that $X_\infty$ is a free $\Lambda$-module of rank 1.
\end{proof}

Therefore, we can choose an isomorphism $i:X_\infty \cong \Lambda$
which induces isomorphisms
$(X_\infty)_{\Gamma_n} \cong \Lambda_n$ 
for each $n$.  Then (\ref{eqn:map3}) becomes
\begin{eqnarray}
\label{eqn:map4}
\hat{E}(\Q_{n,p}) \stackrel{F_n}{\maps} \Lambda_n \maps \Sel_p(E/\Qn)^\dual \maps 0.
\end{eqnarray}
One can verify the commutativity of 
\begin{eqnarray}
\label{diag:1}
\begin{CD}
\hat{E}(\Q_{n,p})  @>F_n>>   \Lambda_n  \\
@V\Tr_{n/n-1}VV      @VV \pin V \\
\hat{E}(\Q_{n-1,p})  @>F_{n-1}>>   \Lambda_{n-1}
\end{CD}
\end{eqnarray}
and
\begin{eqnarray}
\label{diag:2}
\begin{CD}
\hat{E}(\Q_{n,p})  @>F_n>>   \Lambda_n  \\
@Ai_{n-1/n}AA      @AA \nun A \\
\hat{E}(\Q_{n-1,p})  @>F_{n-1}>>   \Lambda_{n-1}
\end{CD}
\end{eqnarray}
where $\Tr_{n/n-1}$ is the trace map,
$\pin$ is the natural projection,
$i_{n-1/n}$ is the natural inclusion and $\nun$ is defined by
$$
\nun(\sigma) = \sum_{\substack{\tau \rightarrow \sigma \\ \tau \in G_n}} \tau
$$
for $\sigma \in G_{n-1}.$  (See \cite[Proposition 6.3]{Iovita} for
a detailed explanation of why these diagrams commute.)

\section{Formal groups}

We now state a result of Kobayashi that describes generators
of $\hat{E}(\Q_{n,p})$ as a Galois module.  

\begin{thm}
\label{thm:formal}
Let $p$ be an odd prime.  
For each $n \geq 0$ there exists $c_n \in \hat{E}(\Q_{n,p})$ such that
\begin{enumerate}
\item $\Tr_{n/n-1} c_n = a_p c_{n-1} - i_{n-2/n-1}(c_{n-2})$ for $n \geq 2$.
\item $\Tr_{1/0} c_1 = \left(  a_p - \frac{p-1}{a_p-2} \right) c_0$.
\end{enumerate}
Furthermore, as a Galois module, $\hat{E}(\Q_{n,p})$ is generated by
$c_n$ and $i_{n-1/n}(c_{n-1})$ for 
$n \geq 1$ and $\hat{E}(\Qp)$ is generated by $c_0$.
\end{thm}

\begin{proof} 
The points $c_n$ were originally constructed by Perrin-Riou in \cite{PR:90}.  
In \cite{Kobayashi}, Kobayashi gives an alternate construction
of these points using Honda theory and 
proves that they generate the formal group as a
Galois module (see \cite[Proposition 8.12]{Kobayashi}).

We point out that Kobayashi assumes that $a_p=0$, but with minor
modifications his arguments would work for any $a_p$ divisible by $p$.
Namely, in the notation of 
\cite{Kobayashi}, one has a formal group $\F :=
\F_{ss}$ whose logarithm is of Honda type $t^2 + p$.  We must replace
$\F$ with a formal group whose logarithm is of Honda type $t^2 - a_p
t + p$.  

Consider the sequence $\{ x_k \}$ defined by
$x_{-1} = 0$, $x_0 = 1$ and 
$$ p x_k - a_p x_{k-1} + x_{k-2} = 0 
$$ for $k
\geq 1$.  Then there exists a formal group $\F(a_p)$ such that 
$$
\log_{\F(a_p)}(X) = \sum_{k=0}^\infty x_k ((X+1)^{p^k} - 1)
$$ 
and its logarithm is of Honda type $t^2 - a_p t + p$ (see \cite[pg. 221]{Honda}).

A second change that needs to made is that Kobayashi chooses an element $\ve
\in p\Zp$ such that $\log_{\F}(\ve)= \frac{p}{p+1}$.  To make the
computations of \cite[Lemma 8.9]{Kobayashi} work out for general $a_p$, we
must choose $\ve \in p\Zp$ such that $\log_{\F(a_p)}(\ve) =
\frac{p}{p+1-a_p}$. With these two modifications, Kobayashi's arguments apply
to this more general setting. 
\end{proof}

\section{$\mu$ and $\lambda$-invariants}

The proof of Theorem \ref{thm:main} will boil
down to understanding the size of certain explicit quotients of $\Lambda_n$.
In this section, we introduce
the notion of $\mu$ and $\lambda$-invariants of elements of
$\Lambda_n$ to help in determining the size of such quotients.

\begin{defn}
For non-zero $f \in \Lambda_n$ the $\mu$-invariant of $f$ is the unique integer
$\mu(f)$ such that $f \in p^{\mu(f)} \Lambda_n - p^{\mu(f)+1} \Lambda_n$.
\end{defn}

Let $I_n$ be the augmentation ideal of $\Lambda_n$ and let
$\widetilde{I}_n$ be the augmentation ideal of $\widetilde{\Lambda}_n :=
\Fp[G_n]$.

\begin{defn} 
For non-zero $f \in \Lambda_n$ the $\lambda$-invariant of $f$
is the unique integer $\lambda(f)$ such that the reduction mod $p$ of
$p^{-\mu(f)} f$ lands in $\widetilde{I}_n^{\lambda(f)} -
\widetilde{I}_n^{\lambda(f)+1}$.  
\end{defn}

\begin{remark}
\label{rmk:cycl}
These $\mu$ and $\lambda$-invariants of elements of $\Lambda_n$ are
related to the standard 
Iwasawa invariants of elements of $\Lambda$.  Namely,
if $f \in \Lambda$ and $f_n$ is its image in $\Lambda_n$
then 
$$
\mu(f) = \mu(f_n) \text{~~and~~}
\lambda(f) = \lambda(f_n) 
$$
if $\lambda(f) < p^n$.  
\end{remark}

Since the ring $\Lambda_n$ is not a domain these invariants do not
share all of the basic properties of standard $\mu$ and
$\lambda$-invariants.  
For instance, since $p \Lambda_n$ is not a prime ideal there exist
$f,g \in \Lambda_n$ such that $\mu(f) = \mu(g) = 0$ but $\mu(f \cdot g) > 0$.
The following simple lemma states some weaker
properties that are true of these invariants.  

\begin{lemma}
\label{lemma:basic}
For $f,g \in \Lambda_n$ we have
\begin{enumerate}
\item $\mu(f \cdot g) \geq \mu(f) + \mu(g)$.
\item 
If $\mu(f \cdot g) = 0$ then $\lambda(f \cdot g) = \lambda(f) +
\lambda(g)$.
\end{enumerate}
 \end{lemma}

These invariants can be used to describe the valuations of elements of
$\Lambda_n$ evaluated at finite order characters as demonstrated in the
following lemma. 

\begin{lemma}
\label{lemma:eval}
Let $f \in \Lambda_n$ and let $\chi$ be a character
of $G_n$ of order $p^n$.  If $\lambda(f) < p^{n-1}(p-1)$ then
$$
\ordp(\chi(f)) = \mu(f) + \frac{\lambda(f)}{p^{n-1}(p-1)}.
$$
\end{lemma}

\begin{proof}
Let $\gamma$ be a generator of $G_n$.  Then $\gamma-1$ is a
generator of the augmentation ideal $I_n$.  From the definitions of
$\mu$ and $\lambda$-invariants, we have that 
$$
f = p^{\mu(f)} \left((\gamma-1)^{\lambda(f)} \cdot u + p \cdot g \right)
$$
for $u \in \Lambda_n^\times$ and $g \in \Lambda_n$.  Hence
\begin{align*}
\ordp(\chi(f)) &= \mu(f) + \min \Big\{  \lambda(f) \cdot \ordp(\chi(\gamma)-1),
			1+\ordp(\chi(g))  \Big\} \\
	&= \mu(f) + \frac{\lambda(f)}{p^{n-1}(p-1)} 
\end{align*}
since $\lambda(f) < p^{n-1}(p-1)$.  
\end{proof}

We will need to understand how these invariants are affected by the
maps $\nun$ and $\pin$.
We first give a lemma that 
describes the relations between these two maps.

\begin{lemma} For $f \in \Lambda_{n-1}$ and $g \in \Lambda_n$ we have
\label{lemma:pinu}
\begin{enumerate}
\item $\pin(\nun(f)) = p \cdot f$
\item 
\label{part:pinutwo}
$\nun(\pin(g)) = \xi_n \cdot g$
\item 
\label{part:pinuthree}
$\im(\nun) = \xi_n \Lambda_n$.
\end{enumerate}
where $\displaystyle \xi_n = \sum_{\sigma^p=1} \sigma \in \Zp[G_n]$.
\end{lemma}

\begin{proof}
%

This lemma follows directly from the definitions.
\end{proof}

We now compute the $\mu$ and $\lambda$-invariant of
the element $\xi_n$ defined in the previous lemma.

\begin{lemma}
\label{lemma:cycl}
We have that $\mu(\xi_n)=0$ and $\lambda(\xi_n) = p^n - p^{n-1}$.
\end{lemma}

\begin{proof} 
Let $\gamma$ be a generator of $G_n$.  Then both
$I_n$ and $\widetilde{I_n}$ are principal generated by $\gamma - 1$. So 
$$
\xi_n = \sum_{\sigma^p=1} \sigma = \sum_{a=0}^{p-1} \gamma^{a p^{n-1}}
= \frac{\gamma^{p^n}-1}{\gamma^{p^{n-1}}-1}
\equiv (\gamma-1)^{p^n-p^{n-1}} \pmod{p}
$$
and hence $\mu(\xi_n)=0$ and $\lambda(\xi_n) = p^n - p^{n-1}$.
\end{proof}

\begin{remark}
If we fix a generator of $G_n$ and thus an isomorphism
$$
\Lambda_n \cong \Zp[[T]] /( (1+T)^{p^n}-1),
$$
the element $\xi_n \in \Lambda_n$ is identified with 
$\Phi_n(1+T)$ where $\Phi_n$ is the $p^n$-th cyclotomic polynomial.  
Note that the computations of the previous lemma agree with the 
computations of the
standard $\mu$ and $\lambda$-invariants of $\Phi_n(1+T)$ 
as predicted by Remark \ref{rmk:cycl}.

\end{remark}

The following proposition summarizes how the Iwasawa
invariants interact with the maps $\nun$ and $\pin$.

\begin{prop} For $f \in \Lambda_{n-1}$ and $g,h \in \Lambda_n$ we have
\label{prop:mulam}
\begin{enumerate}
\item 
\label{part:mulam1}
$\mu(\pin(g)) \geq \mu(g)$ and thus if $\mu(\pin(g)) =0$
then $\mu(g) = 0$.
\item 
\label{part:mulam2}
If $\mu(\pin(g))=\mu(g)$ then $\lambda(\pin(g)) = \lambda(g)$.
\item 
\label{part:mulam3}
$\mu(\nun(f))=\mu(f)$.
\item 
\label{part:mulam4}
$\lambda(\nun(f))=p^{n}-p^{n-1} + \lambda(f)$.
\end{enumerate}
\end{prop}

\begin{proof}
Part \ref{part:mulam1} follows directly from the definitions.
For part \ref{part:mulam2}, 
we have that $\widetilde{g} \in \widetilde{I}_n^a$
if and only if $\pin(\widetilde{g}) \in \widetilde{I}_{n-1}^a$
since these augmentation ideals are principal.  
(Here $\widetilde{g}$ represents the reduction of $g$ mod $p$.)
Thus,
$\lambda(\pin(g)) = \lambda(g)$ since the $\mu$-invariant of 
both of these elements are the same.

For part \ref{part:mulam3}, 
write $f = p^{\mu(f)} f'$ with $\mu(f')=0$.  Then
$\nun(f) = p^{\mu(f)} \nun(f')$ and if we knew that
$\mu(\nun(f'))=0$ then we would have $\mu(\nun(f)) = \mu(f)$.
Hence, we have reduced to the case where $\mu(f) = 0$.  Now pick any  
$g \in \Lambda_n$ such that $\pin(g) = f$.  (Note then by part
\ref{part:mulam1}, $\mu(g) = 0$.)  So
$$
\nun(f)= \nun(\pin(g)) = \xi_n \cdot g
$$ 
by Lemma \ref{lemma:pinu}.\ref{part:pinutwo} and thus
$$
\mu(\nun(f)) = \mu(\xi_n \cdot g) = \mu(g) = 0 = \mu(f).
$$

For the last part, as in part \ref{part:mulam3}, we may assume that
$\mu(f) = 0$.   
Then pick $g \in \Lambda_n$ lifting $f$ and thus
\begin{align*}
\lambda(\nun(f)) &= \lambda(\xi_n \cdot g) &&\\
&= \lambda(\xi_n) + \lambda(g) 
&& \text{(by~part~3~and~Lemma~\ref{lemma:basic})} \\
&= p^n - p^{n-1} + \lambda(\pin(f)) 
&& \text{(by~Lemma~\ref{lemma:cycl})} \\
&= p^n - p^{n-1} + \lambda(f) 
&& \text{(by~part~\ref{part:mulam2})}.
\end{align*}

\end{proof}


We introduce one more lemma which will be useful in the following section.

\begin{lemma}
\label{lemma:domain}
Let $f,g$ be elements of $\Lambda_n$ such that $f \cdot g \in \im(\nun)$.  
If $\mu(f) = 0$ and $\lambda(f) < p^{n-1}$ then $g \in \im(\nun)$.
\end{lemma}

\begin{proof}
By Lemma \ref{lemma:pinu}.\ref{part:pinuthree}, $\im(\nun) = \xi_n 
\Lambda_n$.  Thus, $\im(\nun)$ is a prime ideal in $\Lambda_n$ since
$\Lambda_n / \xi_n \Lambda_n \cong \Zp[\mu_{p^n}]$ which is a domain.
Hence $f \cdot g \in \im(\nun)$ implies that either
$f \in \im(\nun)$ or $g \in \im(\nun)$.  

If 
$f \in \im(\nun)$ then $f = \xi_n h$ for some $h \in \Lambda_n$.
Since $\mu(f) = 0$,
$$
\lambda(f) 
 \geq \lambda(\xi_n) 
 = p^n - p^{n-1} 
 \geq p^{n-1}
$$
by Lemma \ref{lemma:basic}.  This contradicts our hypothesis and
thus $g \in \im(\nun)$.
\end{proof}

\section{Main argument}

Recall the map $F_n : \hat{E}(\Q_{n,p}) \maps \Lambda_n$ defined in
(\ref{eqn:map4}).  For $c_n \in \hat{E}(\Q_{n,p})$ defined in Theorem
\ref{thm:formal}, set 
$$
P_n = F_n(c_n) \in \Lambda_n.
$$ 
The trace relations
between the $c_n$ then yield relations between the $P_n$ by diagrams
(\ref{diag:1}) and (\ref{diag:2}).  We have
\begin{eqnarray}
\label{eqn:tracerel}
\pi_{n+1/n}(P_{n+1}) = a_p P_{n} - \nun(P_{n-1}),
\end{eqnarray}
$$
\pi_{1/0}(P_1) = u P_0 \text{~~with~}u \in \Zpx.
$$
Since $c_n$ and $i_{n-1/n}(c_{n-1})$ generate $\hat{E}(\Q_{n,p})$ as
a Galois module, (\ref{eqn:map4}) yields
\begin{eqnarray}
\label{eqn:explicitsel}
\Lambda_n / (P_n , \nun(P_{n-1}) ) \cong \Sel_p(E/\Qn)^\dual \text{~for~} n \geq 1 \text{~and~}
\end{eqnarray}
$$
\Lambda_0 / (P_0) \cong \Sel_p(E/\Q)^\dual.
$$
Our goal is thus to compute the size of $\Lambda_n/(P_n,\nun(P_{n-1}))$.

We first compute the $\mu$ and $\lambda$-invariants of $P_n$.
For $n \geq 2$, let
$$
q_n =
\begin{cases}
p^{n-1} - p^{n-2} + \cdots + p - 1 &\text{~for~} 2 \mid n \\
p^{n-1} - p^{n-2} + \cdots + p^2 - p &\text{~for~} 2 \nmid n
\end{cases}
$$
and set $q_0 = q_1 = 0$.  

\begin{lemma}
\label{lemma:mulambda}
For $n \geq 0$, 
\begin{enumerate}
\item \label{item:mu} $\mu(P_n) = 0$.
\item \label{item:lambda} $\lambda(P_n) = q_n$.
\end{enumerate}
\end{lemma}

\begin{proof}
We have $\Lambda_0 / (P_0) \cong \Sel_p(E/\Q)^\dual = 0$.  Hence
$P_0$ is a unit and thus $P_1$ is a unit since
$\pi_{1/0}(P_1) = u P_0$ with $u \in \Zpx$,
Therefore, $\mu(P_0) = \mu(P_1) = 0$.
Proceeding by induction, we assume that $\mu(P_k) = 0 $ for $k \leq n$.
We have
\begin{align*}
\mu(\pi_{n+1/n}(P_{n+1})) 
&= \mu(a_p P_n - \nun(P_{n-1})) &&(\text{by~(\ref{eqn:tracerel})})\\
&= \mu(\nun(P_{n-1}))           &&(\text{since~} p \mid a_p)  \\
&= \mu(P_{n-1})         &&(\text{by~Proposition~\ref{prop:mulam}.\ref{part:mulam3}}) \\ 
&= 0.
\end{align*}
Thus, by Proposition \ref{prop:mulam}.\ref{part:mulam1}, $\mu(P_{n+1})=0$
which completes the proof of part \ref{item:mu}.

As for part \ref{item:lambda},
we have already seen that $P_0$ and $P_1$ are units and hence
$\lambda(P_0) = \lambda(P_1) = 0 = q_0 = q_1$.
Again, proceeding by induction, assume that $\lambda(P_k) = q_k$ 
for $k \leq n$.  We have
\begin{align*}
\lambda(\pi_{n+1/n}(P_{n+1})) 
&= \lambda(a_p P_{n} - \nun(P_{n-1}))  \\ 
&= \lambda(\nun(P_{n-1})) \\
&= p^{n} - p^{n-1} + \lambda(P_{n-1})   \hspace{1.5cm}
(\text{by~Proposition~\ref{prop:mulam}.\ref{part:mulam4}}) \\
&= p^{n} - p^{n-1} + q_{n-1} \\
&= q_{n+1}.
\end{align*}
Since we have already seen that
$\mu(\pi_{n+1/n}(P_{n+1})) = 0$,
by Proposition \ref{prop:mulam}.\ref{part:mulam2},
we conclude that $\lambda(P_{n+1}) = \lambda(\pi_{n+1/n}(P_{n+1}))
= q_{n+1}$ 
completing the proof.
\end{proof}

The following lemma will be key in performing the necessary induction 
to compute the size of $\Lambda/(P_n, \nun(P_{n-1}))$.

\begin{lemma}
\label{lemma:exact}
We have an exact sequence
$$
0 \maps \Lambda_{n-1} / J_{n-1} \stackrel{\nun}{\maps}
        \Lambda / J_n \stackrel{\chi}{\longrightarrow} 
	\Zp[\mu_{p^n}] / (\chi(P_n)) \maps 0
$$
where $J_n = (P_n, \nun(P_{n-1}))$ for $n>0$, $J_0 = (P_0)$ 
and $\chi$ is a character of $G_n$ of order $p^n$.
\end{lemma}

\begin{proof}
We check that $\nun(J_{n-1}) \subseteq J_n$ and
that the first map is injective.  The other details are
straightforward to verify.

We have that 
\begin{align*}
\nun(\nu_{n-2/n-1}(P_{n-2}))
&= \nun(a_p P_{n-1} - \pin(P_n)) \hspace{1.0cm} \text{(by~(\ref{eqn:tracerel}))}\\
&= a_p \nun(P_{n-1}) - \xi_n P_n 
\end{align*}
which lies in $J_n$ and thus
$$
\nun(J_{n-1}) = \left(\nun(P_{n-1}), \nun(\nu_{n-2/n-1}(P_{n-2}))\right)
\subseteq J_n.
$$
Thus the first map is well-defined.

To check injectivity, let $f \in \Lambda_{n-1}$ such that
$\nun(f) \in J_n$.  Then
$$
\nun(f) = \alpha \cdot P_n + \beta \cdot \nun(P_{n-1})
$$
and we see that $\alpha \cdot P_n \in \im(\nun)$.  By Lemma \ref{lemma:domain},
$\alpha = \nun(\alpha')$ for some $\alpha' \in \Lambda_{n-1}$
since $\mu(P_n)=0$ and $\lambda(P_n) = q_n < p^{n-1}$.
Hence
$$
\nun(f) = \nun(\alpha') \cdot P_n + \beta \cdot \nun(P_{n-1})
$$
and applying $\pin$ yields
\begin{align*}
p \cdot f &= p \cdot \alpha' \cdot \pin(P_n) + p \cdot \pin(\beta) 
\cdot P_{n-1} \\
    &= p \cdot \alpha' \cdot (a_p P_{n-1} - \nu_{n-2/n-1}(P_{n-2})) 
         + p \cdot \pin(\beta) \cdot P_{n-1}
\end{align*}
which lies in $p J_{n-1}$.  Since $\Lambda_n$ is $p$-torsion free, 
we have that $f \in J_{n-1}$ which establishes
the injectivity of the first map.
\end{proof}

Recall the quantity $e_n$ defined in section \ref{section:intro}.  

\begin{prop}
\label{prop:e_n}
For $n \geq 0$,
$$
\ordp \left(  \# \Lambda_n / J_n \right) = e_n.
$$
\end{prop}

\begin{proof}
For $n=0$ we have $\Lambda_0 / J_0 \cong \Lambda_0 / (P_0)
= 0 = e_0$ since $P_0$ is a unit.
We proceed by induction on $n$.
By direct computation, Lemma \ref{lemma:mulambda} and Lemma \ref{lemma:eval},
we have that
$$
\ord_p\left(\#(\Zp[\mu_{p^n}] / \chi(P_n))\right) = p^{n-1}(p-1) \cdot \ordp(\chi(P_n)) 
= q_n
$$  
where $\chi$ is a character on $G_n$ of order $p^n$.
Therefore, by induction and Lemma \ref{lemma:exact}, we have 
\begin{align*}
\ordp \left(  \# \Lambda_n / J_n \right) &= 
\ordp \left(  \# \Lambda_{n-1} / J_{n-1} \right) +
\ordp\left(\#(\Zp[\mu_{p^n}] / \chi(P_n))\right) \\
&= e_{n-1} + q_n = e_n.
\end{align*}
\end{proof}

\begin{proof}[Proof of Theorem \ref{thm:main}]
By (\ref{eqn:explicitsel}), we have for $n \geq 0$,
$$
\Sel_p(E/\Qn)^\dual \cong \Lambda_n / J_n.
$$  
Hence, by Proposition \ref{prop:e_n},
$$
\ordp(\#\Sel_p(E/\Qn)) = e_n
$$
and, in particular, it is a finite group.  Thus, $E(\Qn)$ is finite
(proving part 1) and $\ordp(\#\Sha(E/\Qn)[p^\infty]) = e_n$  (proving
part 2). 

Now, if $a_p=0$ we have
\begin{align*}
\Tr_{n/m}(c_n) &= \Tr_{n-1/m}(-i_{n-2/n-1}(c_{n-2})) \\
&= -p \Tr_{n-2/m}(c_{n-2}) =
\cdots 
= \pm p^r i_{m-1/m}(c_{m-1})
\end{align*}
for some $r$ when $m$ and $n$ have different parities.  Thus, by diagram (\ref{diag:2}), 
$$\pi_{n/m}(P_n) \in \im(\nu_{m-1/m})$$ and, by Lemma \ref{lemma:pinu}.\ref{part:pinuthree},
$\chi(P_n) = 0$ for $\chi$ of order $p^m$.  Therefore, $P_n \in J_n^\ve$ for $\ve = (-1)^{n+1}$ and
$$
J_n = (P_n , \nun(P_{n-1})) \subseteq J_n^+ + J_n^-.
$$
Then, comparing sizes, we see that
$$
\Sha(E/\Qn)[p^\infty]^\dual \cong \Lambda_n / (P_n, \nun(P_{n-1}))
		    \cong \Lambda_n / (J_n^+ + J_n^-)
$$
completing the proof of part 3.
\end{proof}

\begin{remark}
Note that under Kurihara's hypotheses, \cite[Proposition 1.2]{Kurihara}
implies that
$$
J_n^+ + J_n^- = (\theta_n , \nun(\theta_{n-1}))
$$
where $\theta_n \in \Lambda_n$ is the Mazur-Tate-Teitelbaum element defined
via modular symbols.
Hence part 3 of Theorem \ref{thm:main} is consistent with the isomorphism
$$
\Sha(E/\Qn)[p^\infty]^\dual \cong \Lambda / (\theta_n , \nun(\theta_{n-1}))
$$
proven in \cite{Kurihara}.
\end{remark}

\end{document}